\documentclass[a4paper,11pt]{amsart}
\usepackage[utf8]{inputenc}
\usepackage{amsmath, amssymb, amsthm}
\usepackage{geometry}
\usepackage{enumitem}
\usepackage{booktabs}
\usepackage{hyperref}

\newtheorem{theorem}{Theorem}[section]
\newtheorem{lemma}[theorem]{Lemma}
\newtheorem{proposition}[theorem]{Proposition}

\theoremstyle{definition}

\newcommand{\Q}{\mathbb{Q}}
\newcommand{\Z}{\mathbb{Z}}

\title[On $(P_{n+1}^{(k)})^x-(P_{n-1}^{(k)})^x=2P_m^{(k)}$]{On the exponential Diophantine equation
$\left(P_{n+1}^{(k)}\right)^{x}-\left(P_{n-1}^{(k)}\right)^{x}=2P_{m}^{(k)}$}

\author[F. Faye]{Fallou Faye}
\address{F. Faye, UFR SAT, Universit\'e Alioune Diop, Bambey, 30, S\'en\'egal}
\email{falloufaye1993@gmail.com}

\author[B. Faye]{Bernadette Faye}
\address{B. Faye, UFR SAT, Universit\'e Alioune Diop, Bambey, 30, S\'en\'egal}
\email{bernadette.faye@uadb.edu.sn}


\begin{document}

\begin{abstract}
For the Pell numbers one has $P_{n+1}-P_{n-1}=2P_n$ and $P_{n+1}^2-P_{n-1}^2=2P_{2n}$. We look
for all the solutions of $\big(P_{n+1}^{(k)}\big)^x-\big(P_{n-1}^{(k)}\big)^x=2P_m^{(k)}$ in
positive integers $k\ge2$, $n$, $m$, $x$, where $P^{(k)}$ is the $k$-generalized Pell
sequence. We prove that for $k=2$ the two identities above give all the solutions, and that
for $k\ge3$ the only solutions are $(n,m,x)=(1,1,1),(1,2,2),(2,2,1)$. In contrast with the
analogous equation for $k$-generalized Fibonacci numbers, there are no parametric families of
solutions. The proof uses an auxiliary equation with Fibonacci numbers which does not depend
on $k$, Matveev's lower bound for linear forms in logarithms, the reduction of Dujella and
Peth\H{o}, and, for large $k$, the fact that the dominant root of $P^{(k)}$ is very close to
$\varphi^2$.
\end{abstract}

\maketitle

\section{Introduction}

The Fibonacci numbers satisfy $F_{n+1}^2-F_{n-1}^2=F_{2n}$. Patel and Chaves \cite{PC2021}
asked when $F_{n+1}^x-F_{n-1}^x$ is again a Fibonacci number and solved the equation
$F_{n+1}^x-F_{n-1}^x=F_m$ completely. Apart from trivial cases, the only solutions come from
$x=2$, $m=2n$. G\'omez, Luca and G\'omez \cite{GLG2022} treated $F_{n+1}^x-F_{n-1}^x=F_m^y$.
The same question for the $k$-generalized Fibonacci numbers was studied by Bensella, Patel and
Behloul \cite{BPB2022}, and solved completely by G\'omez, G\'omez and Luca \cite{GGL2024}; in
that case there are two infinite parametric families of solutions, with $n\in\{4,5\}$ and $k$
growing exponentially with the parameter. These families come from the fact that the
$k$-generalized Fibonacci numbers $F_r^{(k)}$ are powers of $2$ for $2\le r\le k+1$.

in this paper, we are  interested in the Pell analogue because of the two identities
\begin{equation}\label{eq:pellid}
P_{n+1}-P_{n-1}=2P_n,\qquad P_{n+1}^2-P_{n-1}^2=(P_{n+1}-P_{n-1})(P_{n+1}+P_{n-1})=2P_{2n},
\end{equation}
where $P_n$ is the $n$-th Pell number (the second one follows from the first and from
$P_{n+1}+P_{n-1}=Q_n$, $P_nQ_n=P_{2n}$, with $Q_n$ the Pell--Lucas numbers). So, for Pell
numbers, the natural question is when $P_{n+1}^x-P_{n-1}^x$ is twice a Pell number. Let $k\ge2$
and let $P^{(k)}=(P_n^{(k)})_n$ be the $k$-generalized Pell sequence.
We study the equation
\begin{equation}\label{eq:main}
\left(P_{n+1}^{(k)}\right)^{x}-\left(P_{n-1}^{(k)}\right)^{x}=2P_{m}^{(k)}
\end{equation}
in positive integers $k\ge2$, $n$, $m$ and $x$. In a companion paper \cite{FF2026s} we treat, in
the same spirit, $\big(P_n^{(k)}\big)^x+\big(P_r^{(k)}\big)^x=P_m^{(k)}$, which generalizes the
Pell identity $P_n^2+P_{n+1}^2=P_{2n+1}$ to arbitrary indices.

In a first version of this work we tried to follow the method of \cite{GGL2024} step by step,
and we obtained two parametric families similar to the families for $k$-Fibonacci numbers. A
numerical check showed that they are not solutions. For example, the first member
$(k,n,m,x)=(7,4,11,3)$ gives $$34^3-5^3=39179 \quad \text{while}\quad 2P_{11}^{(7)}=21854.$$ The reason is that
the property used in \cite{GGL2024} has no analogue here. By a result of K{\i}l{\i}\c{c}
\cite{Kilic2008}, the first $k+1$ non-zero terms of $P^{(k)}$ are not powers of $2$ but their are equal to the 
Fibonacci numbers with odd index, $P_r^{(k)}=F_{2r-1}$ for $1\le r\le k+1$. Hence, when $n\le
k$ and $m\le k+1$, \eqref{eq:main} becomes an equation which does not depend on $k$, namely
$$F_{2n+1}^x-F_{2n-3}^x=2F_{2m-1}$$ and this equation has no solution (see Theorem~\ref{thm:fib}).
This is the main reason why the final answer is so short. In fact, in the range where the
$k$-Fibonacci numbers produce whole families of solutions, the $k$-Pell numbers produce none.

\begin{theorem}\label{thm:main}
Let $k\ge2$, $n\ge1$, $m\ge1$ and $x\ge1$ be integers satisfying \eqref{eq:main}.
\begin{enumerate}[label=\rm(\roman*)]
\item If $k=2$, then $(m,x)=(n,1)$ or $(m,x)=(2n,2)$; conversely, these are solutions for every
$n\ge1$.
\item If $k\ge3$, then $(n,m,x)\in\{(1,1,1),\,(1,2,2),\,(2,2,1)\}$, and these three triples are
solutions for every $k\ge3$.
\end{enumerate}
\end{theorem}

\begin{theorem}\label{thm:fib}
The equation
\begin{equation}\label{eq:fib}
F_{2n+1}^x-F_{2n-3}^x=2F_{2m-1}
\end{equation}
has no solution in integers $n\ge2$, $m\ge1$ and $x\ge2$.
\end{theorem}

We use the following strategy for the proof. The elementary cases ($x=1$, $n=1$, and $k=2$ with $x=2$)
are settled in Section~\ref{sec:elem}, and Theorem~\ref{thm:fib} is proved in
Section~\ref{sec:fib}, with three linear forms in logarithms in $\Q(\sqrt5)$. For an arbitrary $k$
we use the same three forms, now in $\Q(\alpha)$, where $\alpha$ is the dominant root of
$P^{(k)}$. This bounds $x$ in terms of $k$ (Section~\ref{sec:generic}).  This gives first
$k<1.05\times10^{45}$, and six rounds of reduction bring this down to $k\le257$
(Section~\ref{sec:largek}). The remaining values $2\le k\le257$ are treated one by one in
Section~\ref{sec:smallk}. 



\section{Preliminaries}\label{sec:prelim}

\subsection{The generalized Pell sequence}

Let $k \geq 2$ be an integer. The $k$-generalized Pell sequence $P^{(k)}:=(P^{(k)}_n)_{n\geq
2-k}$ is defined by
$$P^{(k)}_n=2P^{(k)}_{n-1}+P^{(k)}_{n-2}+\cdots + P^{(k)}_{n-k}\quad\hbox{for all $n\geq 2$},$$
with the initial conditions $P^{(k)}_{-(k-2)}=P^{(k)}_{-(k-3)}=\cdots=P^{(k)}_{-1}=P^{(k)}_0=0$
and $P^{(k)}_1=1$. For $k=2$ we obtain the Pell sequence $0,1,2,5,12,29,70,\ldots$. The terms of
this sequence are called $k$-generalized Pell numbers. When $k$ is fixed we simply write
$P_n:=P_n^{(k)}$.

K{\i}l{\i}\c{c} \cite{Kilic2008} showed that the first $k+1$ non-zero terms of $P^{(k)}$ are
Fibonacci numbers with odd index, namely
\begin{equation}\label{eq:kilic}
P_n^{(k)}=F_{2n-1} \quad\hbox{for all $1\leq n\leq k+1$}.
\end{equation}
The characteristic polynomial of $P^{(k)}$,
$$
\Psi_k(x) := x^k - 2x^{k-1} - x^{k-2}-\cdots - x - 1,
$$
is irreducible over $\Q$ and has just one root outside the unit circle \cite{BHL2021}. Let
$\alpha := \alpha(k)$ denote that root. It was proved in \cite{BHL2021} that
\begin{equation}\label{eq:alpha}
\varphi^2\left(1-\varphi^{-k}\right)<\alpha<\varphi^2,
\end{equation}
where $\varphi=(1+\sqrt5)/2$ is the golden section. We write $\alpha_1:=\alpha,\alpha_2,\ldots,
\alpha_k$ for the roots of $\Psi_k$; thus $|\alpha_i|<1$ for $2\le i\le k$. A direct computation
gives
\begin{equation}\label{eq:identity}
(x-1)\Psi_k(x)=x^{k+1}-3x^k+x^{k-1}+1,\qquad\hbox{so that}\qquad \alpha^2-3\alpha+1=-\alpha^{1-k}.
\end{equation}
For $x_0=1+\sqrt2$ we have $x_0^2-3x_0+1=1-\sqrt2$, so \eqref{eq:identity} gives
$(x_0-1)\Psi_k(x_0)=1-(\sqrt2-1)x_0^{k-1}=1-x_0^{k-2}\le0$. Since $\Psi_k$ has exactly one
positive root, we get $\alpha\ge1+\sqrt2$, with equality for $k=2$.

For an integer $k\ge2$ we consider the function
\begin{equation}\label{eq:g}
g_{k}(z) = \dfrac{z-1}{(k+1)z^2-3kz+k-1}=\frac{z-1}{k(z^2-3z+1)+z^2-1} \qquad (z \in
\mathbb{C}),
\end{equation}
and we put $g:=g_k(\alpha)$.

\begin{lemma}[{Bravo, Herrera and Luca \cite{BHL2021}}]\label{lem:BHL}
Let $k\geq 2$. Then:
\begin{enumerate}[label=\rm(\alph*)]
\item $\alpha^{n-2}\leq P_n^{(k)}\leq \alpha^{n-1}$ for all $n\geq 1$;
\item $P_n^{(k)} = \sum_{i=1}^{k}g_k(\alpha_i)\alpha_i^n$ for all $n\ge 2-k$;
\item $e_n:=P_n^{(k)} - g\,\alpha^{n}$ satisfies $|e_n|<1/2$ for all $n \geq 2 - k$;
\item $0.276< g<0.5$.
\end{enumerate}
\end{lemma}

\begin{lemma}[{Bravo and Herrera \cite{BH2020}; see also \cite{BHL2021b}}]\label{lem:gconj}
For every $k\ge2$ and $2\leq i\leq k$ we have $|g_{k}(\alpha_i)|<1$.
\end{lemma}

We shall also use the following elementary facts. For $n\ge2$ we have $P_n\ge2P_{n-1}$, since
$P_n=2P_{n-1}+(\hbox{non-negative terms})$; hence, for $n\ge2$,
\begin{equation}\label{eq:five}
P_{n+1}\ge2P_n+P_{n-1}\ge5P_{n-1}.
\end{equation}
The sequence $(P_n)_{n\ge1}$ is strictly increasing, with $P_1=1$, $P_2=2$ and $P_3=5$ for
every $k\ge2$.

\begin{lemma}\label{lem:pow2}
Let $k\ge2$ and $n\ge1$. If $P_n^{(k)}=2^a$ with $a\ge1$, then $n=2$ and $a=1$.
\end{lemma}
\begin{proof}
See Lemma 2.3 in \cite{FF2026s}.
\end{proof}

\subsection{Approximation for large $k$}

For large $k$ the dominant root $\alpha$ is very close to $\varphi^2$, and $P_n^{(k)}$ remains very
close to $F_{2n-1}$ (in relative terms) also for $n>k+1$, as long as $n$ is small compared with
$\varphi^{2k}$. We put
$$
c:=\frac1{\varphi^2+1}=\frac1{\sqrt5\,\varphi}=\frac{5-\sqrt5}{10}=0.2763\ldots
$$
For instance $F_{2n-1}=c\varphi^{2n}+\varphi^{-2n+1}/\sqrt5$, since $c\varphi=1/\sqrt5$.

\begin{lemma}\label{lem:approx}
Let $k\ge30$ and write $\theta:=\varphi^2-\alpha$, $\delta:=\theta/\varphi^2$ and $g=c(1+\eta)$.
Then
$$
0<\theta<1.1711\,\varphi^{-2k},\qquad 0<\delta<0.4474\,\varphi^{-2k},\qquad
|\eta|\le(0.45k+1.07)\,\varphi^{-2k}.
$$
Consequently, for every $j\ge0$,
\begin{equation}\label{eq:zeta}
g\alpha^j=c\,\varphi^{2j}(1+\zeta_j)\qquad\hbox{with}\qquad
|\zeta_j|\le\big(0.45(k+j)+1.07\big)\varphi^{-2k}.
\end{equation}
\end{lemma}
\begin{proof}
See the proof of Lemma 2.5 in \cite{FF2026s}.
\end{proof}


\subsection{Notations and terminology from algebraic number theory}

Let $\eta$ be an algebraic number of degree $d$ with minimal primitive polynomial over the
integers
$$
a_0x^{d}+ a_1x^{d-1}+\cdots+a_d = a_0\prod_{i=1}^{d}(x-\eta^{(i)}),
$$
where the leading coefficient $a_0$ is positive and the $\eta^{(i)}$'s are the conjugates of
$\eta$. Then the \textit{logarithmic height} of $\eta$ is given by
$$
h(\eta):=\dfrac{1}{d}\left( \log a_0 + \sum_{i=1}^{d}\log\left(\max\{|\eta^{(i)}|,
1\}\right)\right).
$$
In particular, if $\eta=p/q$ is a rational number with $\gcd (p,q)=1$ and $q>0$, then
$h(\eta)=\log\max\{|p|, q\}$. The following properties of $h(\cdot)$ will be used without
reference:
\begin{eqnarray*}
h(\eta\pm \gamma) &\leq& h(\eta) +h(\gamma) +\log 2,\\
h(\eta\gamma^{\pm 1})&\leq & h(\eta) + h(\gamma),\\
h(\eta^{s}) &=& |s|h(\eta) \qquad (s\in\mathbb{Z}).
\end{eqnarray*}
Since $\Psi_k$ is monic, irreducible, and has exactly one root of modulus greater than $1$, we
have $h(\alpha)=(\log\alpha)/k$. For $g$ we use the bound of Bravo, Herrera and Luca
\cite[Lemmas~3 and 6]{BHL2021b}:
\begin{equation}\label{eq:hg}
h(g)<H(k):=4\log\varphi+\log(k+1)\qquad(k\ge2).
\end{equation}
In $\Q(\sqrt5)$ we shall use $h(2)=\log2$, $h(\sqrt5)=\frac12\log5$, $h(\varphi)=\frac12\log\varphi$,
$h(2/\sqrt5)=\frac12\log5$ (minimal polynomial $5y^2-4$), and
$h(2c)=\frac12\log\big(5(1+1/\sqrt5)\big)=0.98953\ldots$, since $2c=1-1/\sqrt5$ has minimal
polynomial $5y^2-10y+4$.

\subsection{Linear forms in logarithms and the reduction lemma}

We use the following result of Matveev \cite{Matveev2000}.

\begin{theorem}[Matveev]\label{thm:matveev}
Let $\gamma_1,\ldots,\gamma_t$ be positive real algebraic numbers in a real algebraic number
field $\mathbb{K}$ of degree $D$, $b_1,\ldots,b_t$ be nonzero integers, and assume that
$$
\Lambda:=\gamma_1^{b_1}\cdots\gamma_t^{b_t} - 1
$$
is nonzero. Then
$$
\log |\Lambda| > -1.4\times 30^{t+3}\times t^{4.5}\times D^{2}(1+\log D)(1+\log B)A_1\cdots A_t,
$$
where $B\geq\max\{|b_1|, \ldots, |b_t|\}$ and $A_i \geq \max\{Dh(\gamma_i), |\log\gamma_i|,
0.16\}$ for all $i=1,\ldots,t$.
\end{theorem}

We write $C(D,t):=1.4\times30^{t+3}\times t^{4.5}\times D^2(1+\log D)$.

For a real number $r$ we denote by $\| r\|$ the distance from $r$ to the nearest integer. The
following is the reduction lemma of Dujella and Peth\H{o} \cite[Lemma~5(a)]{DP1998}; we state it
with a real exponent $w$, the proof being the same.

\begin{lemma}[Dujella and Peth\H{o}]\label{lem:dp}
Let $M$ be a positive integer, let $\tau$ be an irrational number and let $p/q$ be a convergent
of the continued fraction of $\tau$ with $q>6M$. Let $A,B,\mu$ be real numbers with $A>0$ and
$B>1$, and put $\varepsilon:=\|\mu q\|-M\|\tau q\|$. If $\varepsilon>0$, then there is no
solution of
$$
0<|u\tau-v+\mu|<AB^{-w}
$$
in integers $u,v$ and a real number $w$ with $1\le u\le M$ and $w\ge\log(Aq/\varepsilon)/\log B$.
\end{lemma}

In practice we take the first convergent with $q>6M$, and the next ones if $\varepsilon\le0$; we
never needed more than one extra convergent. We shall also use, without further mention, that
if $|\Lambda|<1/2$ and $\Lambda=e^{\Gamma}-1$ with $\Gamma$ real, then $|\Gamma|<2|\Lambda|$. When
$|\Lambda|\ge1/2$, the upper bound for $|\Lambda|$ that we use already gives a bound for the
exponent $w$ which is smaller than the bounds obtained with Lemma~\ref{lem:dp}.

\section{Elementary cases}\label{sec:elem}

Let $(k,n,m,x)$ be a solution of \eqref{eq:main} in positive integers with $k\ge2$.

\begin{lemma}[The case $x=1$]\label{lem:x1}
If $x=1$, then $m=n$ when $k=2$, and $(n,m)\in\{(1,1),(2,2)\}$ when $k\ge3$.
\end{lemma}
\begin{proof}
For $k=2$, $P_{n+1}-P_{n-1}=2P_n$ by \eqref{eq:pellid}, so $P_m=P_n$ and $m=n$. Let $k\ge3$. For
$n=1$, $P_2-P_0=2=2P_1$, and for $n=2$, $P_3-P_1=4=2P_2$. Let $n\ge3$. Then
$$
P_{n+1}-P_{n-1}=2P_n+P_{n-2}+P_{n-3}+\cdots+P_{n+1-k}\ge2P_n+P_{n-2}>2P_n,
$$
since $k\ge3$ and $P_{n-2}\ge P_1=1$, and $P_{n+1}-P_{n-1}<2P_{n+1}$. Since $(P_r)_{r\ge1}$ is
strictly increasing, $2P_m$ cannot lie strictly between $2P_n$ and $2P_{n+1}$.
\end{proof}

\begin{lemma}[The case $n=1$]\label{lem:n1}
If $n=1$, then $(m,x)\in\{(1,1),(2,2)\}$.
\end{lemma}
\begin{proof}
Since $P_2=2$ and $P_0=0$, \eqref{eq:main} becomes $P_m=2^{x-1}$. For $x=1$ we get $m=1$; for
$x\ge2$, Lemma~\ref{lem:pow2} gives $m=2$ and $x=2$.
\end{proof}

\begin{lemma}[The case $k=2$, $x=2$]\label{lem:k2x2}
If $k=2$ and $x=2$, then $m=2n$.
\end{lemma}
\begin{proof}
By \eqref{eq:pellid}, $P_{n+1}^2-P_{n-1}^2=2P_{2n}$, so $P_m=P_{2n}$ and $m=2n$.
\end{proof}

The triples listed in Theorem~\ref{thm:main} are the solutions. This follows from
\eqref{eq:pellid} for $k=2$, and from $P_2-P_0=2P_1$, $P_2^2-P_0^2=2P_2$ and $P_3-P_1=2P_2$ for
$k\ge3$. Therefore it remains to show that \eqref{eq:main} has no solution with
\begin{equation}\label{eq:remaining}
n\ge2,\qquad x\ge2,\qquad (k,x)\ne(2,2).
\end{equation}

\begin{lemma}[Size of $m$]\label{lem:size}
If $n\ge2$ and $x\ge2$, then $(n-1)x+1\le m\le nx+1$. In particular $m>x$, and the integer
$$
j:=(n+1)x-m
$$
satisfies $x-1\le j\le2x-1$.
\end{lemma}
\begin{proof}
By Lemma~\ref{lem:BHL}(a), $2\alpha^{m-2}\le2P_m<P_{n+1}^x\le\alpha^{nx}$, so $m-2<nx$ and $m\le
nx+1$. By \eqref{eq:five}, $P_{n-1}^x\le5^{-x}P_{n+1}^x$, hence
$2\alpha^{m-1}\ge2P_m\ge(1-5^{-x})P_{n+1}^x\ge\frac{24}{25}\alpha^{(n-1)x}$. Thus
$$m-1\ge(n-1)x-\log(25/12)/\log\alpha>(n-1)x-1$$ because $\log(25/12)<\log(1+\sqrt2)$.
Hence $m\ge(n-1)x+1$.
\end{proof}

\section{Proof of Theorem~\ref{thm:fib}}\label{sec:fib}

We begin with the case in which all the indices are at most $k+1$, which reduces to
\eqref{eq:fib}. Everything takes place in $\Q(\sqrt5)$, the heights are small, and the argument
shows on a simple example the steps that we shall repeat later for an arbitrary $k$.

In this section $(n,m,x)$ is a solution of \eqref{eq:fib} with $n\ge2$ and $x\ge2$. We use
$F_r=(\varphi^r-(-\varphi)^{-r})/\sqrt5$ and $\varphi^{r-2}\le F_r\le\varphi^{r-1}$ for $r\ge1$.

Since $F_{r+4}=3F_{r+1}+2F_r\ge5F_r$, we have $F_{2n-3}\le F_{2n+1}/5$.
Exactly as in Lemma~\ref{lem:size} (with $\varphi^2$ in place of $\alpha$) we get $(n-1)x+1\le
m\le nx+1$. We put
$$
e:=2m-1-(2n+1)x,\qquad\hbox{so that}\qquad -3x+1\le e\le-x+1 .
$$

	Write $a:=\varphi^{-4n-2}$ and $b:=\varphi^{-4n+6}$, so that
$F_{2n+1}=\frac{\varphi^{2n+1}}{\sqrt5}(1+a)$ and $F_{2n-3}=\frac{\varphi^{2n-3}}{\sqrt5}(1+b)$,
and $2F_{2m-1}=\frac{2}{\sqrt5}\big(\varphi^{2m-1}+\varphi^{-2m+1}\big)$. We consider
\begin{align*}
\Phi_1&:=\frac2{\sqrt5}\,\varphi^{2m-1}F_{2n+1}^{-x}-1,\\
\Phi_2&:=2\cdot\sqrt5^{\,x-1}\varphi^{e}-1,\\
\Phi_3&:=2\cdot\sqrt5^{\,x-1}\varphi^{e}\,(1-\varphi^{-4x})^{-1}-1.
\end{align*}
Dividing \eqref{eq:fib} by $F_{2n+1}^x$, we get
$\Phi_1=-(F_{2n-3}/F_{2n+1})^x-\frac2{\sqrt5}\varphi^{-2m+1}F_{2n+1}^{-x}$, hence
\begin{equation}\label{eq:Phi1}
|\Phi_1|<1.4\cdot5^{-x}.
\end{equation}
Dividing \eqref{eq:fib} by $(\varphi^{2n+1}/\sqrt5)^x$, we get
$$
\Phi_2=\big[(1+a)^x-1\big]-\varphi^{-4x}(1+b)^x-r,\qquad
r:=\frac2{\sqrt5}\varphi^{-2m+1}\Big(\frac{\sqrt5}{\varphi^{2n+1}}\Big)^x .
$$
Assume that
\begin{equation*}\label{eq:starF}
x\,\varphi^{-4n+6}\le0.05 . \tag{$12a$}
\end{equation*}
Then $(1+a)^x-1\le1.026\,xa$, $(1+b)^x\le e^{0.05}$, and, using $m\ge(n-1)x+1$, $r\le0.014\,xa$.
This gives
\begin{equation}\label{eq:Phi2}
|\Phi_2|<2.7\max\{x\varphi^{-4n-2},\ \varphi^{-4x}\}.
\end{equation}
In the same way, the numerator of $\Phi_3$ is
$[(1+a)^x-1]-\varphi^{-4x}[(1+b)^x-1]-r$. Note that $\varphi^{-4x}b\le a$ for $x\ge2$, and the denominator is
$1-\varphi^{-4x}\ge1-\varphi^{-8}$; hence, under \eqref{eq:starF},
\begin{equation}\label{eq:Phi3}
|\Phi_3|<2.66\,x\,\varphi^{-4n-2}.
\end{equation}

The norms from $\Q(\sqrt5)$ to $\Q$ satisfy
$|N(\varphi)|=1$, $|N(\sqrt5)|=5$, $|N(2/\sqrt5)|=4/5$. If $\Phi_1=0$, then
$|N(2\varphi^{2m-1}/\sqrt5)|=4/5$ would be equal to $F_{2n+1}^{2x}\in\Z$, which is impossible. If
$\Phi_2=0$, then $1=|N(2\sqrt5^{\,x-1}\varphi^e)|=4\cdot5^{x-1}$, which is impossible. If
$\Phi_3=0$, then $$4\cdot5^{x-1}=|N(1-\varphi^{-4x})|=|(1-\varphi^{-4x})(1-\varphi^{4x})|=L_{4x}-2
=5F_{2x}^2$$, that is, $F_{2x}^2=4\cdot5^{x-2}$. This is false for $x=2$ ($F_4=3$); for $x\ge3$ it
would give $5\mid F_{2x}$, hence $5\mid x$, and then $2v_5(x)=2v_5(F_{2x})=x-2$, which is
impossible since $x-2>2\log x/\log5$ for $x\ge5$. Hence $\Phi_1\Phi_2\Phi_3\ne0$.

We apply Theorem~\ref{thm:matveev} in $\mathbb K=\Q(\sqrt5)$, $D=2$:
\begin{itemize}
\item to $\Phi_2$, with $t=3$, $(\gamma_1,\gamma_2,\gamma_3)=(2,\sqrt5,\varphi)$,
$(b_1,b_2,b_3)=(1,x-1,e)$, $B=3x$, and $(A_1,A_2,A_3)=(2\log2,\log5,\log\varphi)$. We get
$\log|\Phi_2|>-K_2(1+\log3x)$, with
$$K_2<1.0413\times10^{12}$$
\item to $\Phi_3$, with $t=4$, the same data and $\gamma_4=1-\varphi^{-4x}$, $b_4=-1$,
$A_4=\log(\varphi^{4x}-1)<4x\log\varphi$ where the conjugate of $\gamma_4$ is $1-\varphi^{4x}$. More, 
$$\log|\Phi_3|>-4x\log\varphi\cdot K_3(1+\log3x),$$ $K_3:=C(2,4)<1.140\times10^{14}$
\item to $\Phi_1$, with $t=3$, $(\gamma_1,\gamma_2,\gamma_3)=(2/\sqrt5,\varphi,F_{2n+1})$,
$(b_1,b_2,b_3)=(1,2m-1,-x)$, $B=2m$, $(A_1,A_2,A_3)=(\log5,\log\varphi,2\log F_{2n+1})$ with
$2\log F_{2n+1}\le4n\log\varphi$ where $$\log|\Phi_1|>-4n\log\varphi\cdot K_1(1+\log2m)$$
$K_1:=C(2,3)\log5\log\varphi<7.511\times10^{11}$.
\end{itemize}

Suppose first that \eqref{eq:starF} holds and that
$\varphi^{-4x}\ge x\varphi^{-4n-2}$. Then \eqref{eq:Phi2} and the bound for $\Phi_2$ give
$4x\log\varphi<\log2.7+K_2(1+\log3x)$, hence $x<X_0:=1.7634\times10^{13}$; then \eqref{eq:Phi3}
and the bound for $\Phi_3$ give $n<6.553\times10^{28}$. Otherwise, either \eqref{eq:starF} fails,
and then $4n<6+\log(20x)/\log\varphi$, or $x\varphi^{-4n-2}>\varphi^{-4x}$, and then
$(4n+2)\log\varphi<\log(2.7x)+K_2(1+\log3x)$. In both cases $n$ is bounded by an explicit function
of $x$; inserting this and $m\le nx+1$ in \eqref{eq:Phi1} and the bound for
$\Phi_1$, we get $x<3.0554\times10^{27}$. Therefore, in all cases,
\begin{equation}\label{eq:fibabs}
x<3.06\times10^{27}.
\end{equation}

For the reduction phase, we put $\Gamma_2:=\log2+(x-1)\log\sqrt5+e\log\varphi$, so that
$\Phi_2=e^{\Gamma_2}-1$. With $\tau:=\log\sqrt5/\log\varphi$ and $\mu:=\log2/\log\varphi$,
\eqref{eq:Phi2} gives (if $|\Phi_2|<1/2$)
$$
0<|(x-1)\tau-(-e)+\mu|<\frac{5.4}{\log\varphi}\,e^{-w},\qquad
w:=\min\{4x\log\varphi,\ (4n+2)\log\varphi-\log x\}.
$$
Lemma~\ref{lem:dp} with $u=x-1\le M$ gives $w<W$ for an explicit $W$. If $w=4x\log\varphi$, then
$x<W/(4\log\varphi)$, and for each such $x$ and each $e\in[-3x+1,-x+1]$ the number $\Phi_3$ is
explicit and non-zero, so that \eqref{eq:Phi3} bounds $n$. If $w$ is the other term, or if
\eqref{eq:starF} fails, $n$ is bounded directly. Once $n\le N$, we apply Lemma~\ref{lem:dp} to
$\Gamma_1/\log F_{2n+1}$, where $$\Gamma_1=\log(2/\sqrt5)+(2m-1)\log\varphi-x\log F_{2n+1}$$ with
$u=2m-1\le2(NM+1)$, $A=2.8/\log F_{2n+1}$, $B=5$, separately for each $2\le n\le N$; this bounds
$x$. Table~\ref{tab:fibred} shows the three passes we made.

\begin{table}[ht]
\centering
\caption{Reduction for equation \eqref{eq:fib}. Here $M$ is the bound for $x$ used in the pass,
$W$ comes from $\Phi_2$, $N$ is the resulting bound for $n$ and the last column is the bound for
$x$ given by $\Phi_1$.}\label{tab:fibred}
\begin{tabular}{@{}rrrrrrr@{}}
\toprule
pass & $M$ & $W$ & $x$ if $w=4x\log\varphi$ & $n$ from $\Phi_3$ & $N$ & $x\le$ \\
\midrule
1 & $3.06\times10^{27}$ & $69.18$ & $35$ & $3.8$ & $68$ & $46$\\
2 & $47$ & $11.00$ & $5$ & $2.3$ & $7$ & $7$\\
3 & $8$ & $7.42$ & $3$ & $1.8$ & $4$ & $7$\\
\bottomrule
\end{tabular}
\end{table}

As usual, almost all the work is done by the first pass, which brings the bound for $x$ from
$3\times10^{27}$ down to $46$. Thus $n\le4$ and $x\le7$. A direct search with $2\le n\le60$ and $2\le x\le60$ shows that
\eqref{eq:fib} has no solution. This proves Theorem~\ref{thm:fib}.\qed

\section{Linear forms in $\Q(\alpha)$ and a bound for $x$}\label{sec:generic}

The forms below are those of Section~\ref{sec:fib}, with $g$
and $\alpha$ in the place of $1/\sqrt5$ and $\varphi^2$. We then obtained that we work is now in a
field of degree $k$ and that the height of $g$ grows like $\log k$, so that the bound for $x$ will
depend on $k$.

In this section $k\ge2$ and $(k,n,m,x)$ is a solution of \eqref{eq:main} satisfying
\eqref{eq:remaining}. We use the notation of Lemma~\ref{lem:BHL}, and $j=(n+1)x-m$ is as in
Lemma~\ref{lem:size}. We consider
\begin{align*}
\Lambda_1&:=2g\,\alpha^{m}P_{n+1}^{-x}-1,\\
\Lambda_2&:=2g^{1-x}\alpha^{-j}-1,\\
\Lambda_3&:=2g^{1-x}\alpha^{-j}(1-\alpha^{-2x})^{-1}-1,
\end{align*}
and the condition
\begin{equation*}\label{eq:star}
36.4\,x\le\alpha^{n-1}. \tag{$15a$}
\end{equation*}

\begin{lemma}\label{lem:Lambda}
\begin{enumerate}[label=\rm(\alph*)]
\item $0<|\Lambda_1|<2\cdot5^{-x}$.
\item If \eqref{eq:star} holds, then $0<|\Lambda_2|<2.6\max\{\alpha^{-2x},\ x\alpha^{-n}\}$ and
$|\Lambda_3|<1.5\,x\,\alpha^{-n}$.
\end{enumerate}
\end{lemma}
\begin{proof}
(a) Since $2g\alpha^m=2P_m-2e_m$, dividing \eqref{eq:main} by $P_{n+1}^x$ gives
$$\Lambda_1=-(P_{n-1}/P_{n+1})^x-2e_mP_{n+1}^{-x}.$$ By \eqref{eq:five} and $P_{n+1}\ge5$ we get
$|\Lambda_1|<2\cdot5^{-x}$. If $\Lambda_1=0$, then $2g\alpha^m=P_{n+1}^x$. Applying an
automorphism of the Galois closure of $\Q(\alpha)$ sending $\alpha$ to $\alpha_2$ and using
Lemma~\ref{lem:gconj}, we get $$P_{n+1}^x=2|g_k(\alpha_2)||\alpha_2|^m<2,$$ a contradiction.

(b) Put $U:=g\alpha^{n+1}$ and $V:=g\alpha^{n-1}$, and write $P_{n+1}=U(1+\varepsilon_1)$,
$P_{n-1}=V(1+\varepsilon_2)$. By Lemma~\ref{lem:BHL}(c),(d),
$|\varepsilon_1|<1.812\,\alpha^{-n-1}$ and $|\varepsilon_2|<1.812\,\alpha^{-n+1}$, so that, by
\eqref{eq:star}, $x|\varepsilon_2|\le0.05$; moreover \eqref{eq:star} implies $\alpha^{n-1}\ge72.8$,
hence $n\ge6$ and $V>20$. Dividing \eqref{eq:main} by $U^x$ and using $V/U=\alpha^{-2}$,
$$
\Lambda_2=\big[(1+\varepsilon_1)^x-1\big]-\alpha^{-2x}(1+\varepsilon_2)^x-2e_mU^{-x}.
$$
Now $$|(1+\varepsilon_1)^x-1|\le e^{x|\varepsilon_1|}-1\le1.026\,x|\varepsilon_1|<0.78\,x\alpha^{-n}$$
as $\alpha\ge1+\sqrt2$, $(1+|\varepsilon_2|)^x<1.052$, and
$U^{-x}=\alpha^{-2x}V^{-x}<0.0025\,\alpha^{-2x}$. This gives the bound for $\Lambda_2$. In the
same way the numerator of $\Lambda_3$ is
$$[(1+\varepsilon_1)^x-1]-\alpha^{-2x}[(1+\varepsilon_2)^x-1]-2e_mU^{-x},$$ of absolute value less
than $1.026\,x(|\varepsilon_1|+\alpha^{-2x}|\varepsilon_2|)+U^{-x}<0.91\,x\alpha^{-n}$, and its
denominator is $1-\alpha^{-2x}\ge1-\alpha^{-4}>0.97$; hence $|\Lambda_3|<1.5x\alpha^{-n}$.
Finally, if $\Lambda_2=0$, then $g^{x-1}\alpha^j=2$; conjugating as in (a) we get
$2=|g_k(\alpha_2)|^{x-1}|\alpha_2|^j<1$, since $j\ge x-1\ge1$, a contradiction.
\end{proof}

\begin{lemma}\label{lem:zero}
\begin{enumerate}[label=\rm(\alph*)]
\item For $k=2$, $x\ge2$ and $x-1\le j\le2x-1$, we have $\Lambda_3=0$ if and only if $(x,j)=(2,2)$.
\item For $2\le k\le257$, $2\le x\le64$, $x-1\le j\le2x-1$ and $(k,x,j)\ne(2,2,2)$, we have
$|\Lambda_3|>1.22\times10^{-3}$.
\end{enumerate}
\end{lemma}
\begin{proof}
(a) For $k=2$ we have $\alpha=1+\sqrt2$ and $g=1/(2\sqrt2)$. For $(x,j)=(2,2)$,
$$2g^{-1}\alpha^{-2}=4\sqrt2\,\alpha^{-2}=1-\alpha^{-4},$$ because
$\alpha^2-\alpha^{-2}=(\alpha-\alpha^{-1})(\alpha+\alpha^{-1})=2\cdot2\sqrt2$. This correspond to  the identity
$$P_{n+1}^2-P_{n-1}^2=2P_{2n}$$ seen in $\Q(\sqrt2)$. Conversely, if $\Lambda_3=0$, taking norms from
$\Q(\sqrt2)$ to $\Q$ gives $2^{3x-1}=|N(1-\alpha^{-2x})|=Q_{2x}-2$, where
$Q_r=\alpha^r+(-\alpha)^{-r}$ are the Pell--Lucas numbers. For $x$ even, $Q_{2x}-2=8P_x^2$, so
$P_x$ is a power of $2$. Thus $x=2$ by Lemma~\ref{lem:pow2}, and then $j=2$. For $x$ odd,
$Q_{2x}-2=Q_x^2$, and $Q_x$ is twice an odd number larger than $1$ for $x\ge3$, which is
impossible.

(b) This is a direct computation (Appendix~\ref{app:comp}). The minimum is attained for $k=4$,
$x=62$, $j=80$, where $|\Lambda_3|>1.22\times10^{-3}$.
\end{proof}

Part (a) explains the only zero we have got. Since $(k,x)=(2,2)$ is excluded by
\eqref{eq:remaining}, it causes no problem. Part (b) covers all the pairs $(x,j)$ for which
$\Lambda_3$ is used in Section~\ref{sec:smallk}, since there $x\le64$.

Now we can apply Matveev for $\Lambda_1$ and $\Lambda_2$. Let $\mathbb K=\Q(\alpha)$, $D=k$. For
$\Lambda_1$ we take $t=3$, $(\gamma_1,\gamma_2,\gamma_3)=(2g,\alpha,P_{n+1})$,
$(b_1,b_2,b_3)=(1,m,-x)$ and $B=m$ (since $m>x$). By \eqref{eq:hg} we may take
$A_1=k(\log2+H(k))$, $A_2=\log\alpha$ and $A_3=k\log P_{n+1}\le kn\log\alpha$. For $\Lambda_2$
we take $(\gamma_1,\gamma_2,\gamma_3)=(2,g,\alpha)$, $(b_1,b_2,b_3)=(1,1-x,-j)$, $B=2x$,
$A_1=k(\log2+H(k))$, $A_2=kH(k)$ and $A_3=\log\alpha$. Using
$\log\alpha<2\log\varphi$, we get
\begin{align}
x\log5-\log2&<E_1(k)\,n\,(1+\log m),\label{eq:M1}\\
\min\{2x\log\alpha,\ n\log\alpha-\log x\}-\log2.6&<E_2(k)(1+\log2x),\label{eq:M2}
\end{align}
where
$$
E_1(k):=C(k,3)\,k^2(\log2+H(k))\,(2\log\varphi)^2,\qquad
E_2(k):=C(k,3)\,k^2(\log2+H(k))\,H(k)\,2\log\varphi,
$$
the second one under \eqref{eq:star}. 

\begin{proposition}\label{prop:X}
Let $X(k)$ be the supremum of the real numbers $x\ge2$ such that either
$$2x\log(1+\sqrt2)<\log2.6+E_2(k)(1+\log2x),$$ or
$$x\log5<\log2+E_1(k)N_s(x)(1+\log(N_s(x)x+1)),$$ where
$$
N_s(x):=1+\frac{\log(36.4x)+\log2.6+E_2(k)(1+\log2x)}{\log(1+\sqrt2)}.
$$
Then every solution of \eqref{eq:main} satisfying \eqref{eq:remaining} has $x<X(k)$. 
\end{proposition}
\begin{proof}
If \eqref{eq:star} holds and the minimum in \eqref{eq:M2} is $2x\log\alpha$, then $x$ satisfies
the first inequality. Otherwise, either \eqref{eq:star} fails or the minimum is
$n\log\alpha-\log x$; in both cases $n<N_s(x)$. Then $m\le nx+1$ by Lemma~\ref{lem:size}, and
\eqref{eq:M1} gives the second inequality.
\end{proof}

\section{Large values of $k$}\label{sec:largek}

 The idea of this
section is to exploit the error made by replacing $\alpha$ by $\varphi^2$. This is of size about
$k\varphi^{-2k}$, which is negligible as soon as $k$ is large, and the forms we
obtain no longer depend on $k$.

In this section $k\ge30$. Lemma~\ref{lem:approx} allows us to replace $g$ by $c=1/(\sqrt5\varphi)$
and $\alpha$ by $\varphi^2$ in $\Lambda_1,\Lambda_2,\Lambda_3$. The new forms live in
$\Q(\sqrt5)$, and their coefficients do not depend on $k$. Since $c^{1-x}=(\sqrt5\varphi)^{x-1}$,
the forms are
\begin{align*}
L_1&:=2c\,\varphi^{2m}F_{2n+1}^{-x}-1,\\
L_2&:=2\sqrt5^{\,x-1}\varphi^{x-1-2j}-1,\\
L_3&:=2\sqrt5^{\,x-1}\varphi^{x-1-2j}(1-\varphi^{-4x})^{-1}-1 .
\end{align*}
Note that $L_2$ and $L_3$ are the forms $\Phi_2$ and $\Phi_3$ of Section~\ref{sec:fib}, with
$e=x-1-2j$. In particular they are non-zero, and Matveev gives the same lower bounds as for
$\Phi_2$, $\Phi_3$. Note that here $|x-1-2j|\le3x$. The form $L_1$ is non-zero because
$|N(2c\varphi^{2m})|=4/5$ is not an integer. Put
$$
E_k(x):=1.011\,x\,(0.45k+2.02)\,\varphi^{-2k}.
$$

\begin{lemma}\label{lem:transfer}
Let $k\ge30$, and let $(k,n,m,x)$ be a solution of \eqref{eq:main} which satisfies
\eqref{eq:remaining} and $2.022\,E_k(x)\le0.01$. Put $\lambda_*:=-\log(2.022E_k(x))$.
\begin{enumerate}[label=\rm(\alph*)]
\item If \eqref{eq:star} holds, then $|L_2|<e^{-\lambda_2}$ and $|L_3|<e^{-\lambda_3}$, where
\begin{align*}
\lambda_2&:=\min\{2x\log\alpha-\log5.26,\ n\log\alpha-\log(5.26x),\ \lambda_*\},\\
\lambda_3&:=\min\{n\log\alpha-\log(3.04x),\ \lambda_*\}.
\end{align*}
\item If $n\le k$ and $|\zeta_m|\le0.01$, then $|L_1|<e^{-\lambda_1}$, where
$$
\lambda_1:=\min\big\{x\log5-\log4.1,\ -\log\big(2.1(0.45(k+m)+1.07)\varphi^{-2k}\big)\big\}.
$$
\end{enumerate}
\end{lemma}
\begin{proof}
(a) By Lemma~\ref{lem:approx}, $2g^{1-x}\alpha^{-j}=(L_2+1)\rho$ with
$\rho=(1+\eta)^{1-x}(1-\delta)^{-j}$, and $$|\log\rho|\le1.0001\big((x-1)|\eta|+j\delta\big)
\le1.0001\,x(0.45k+1.965)\varphi^{-2k}.$$ We recall $j\le2x-1$ and $|\rho-1|\le E_k(x)$. Thus
$$|L_2|=|\Lambda_2+1-\rho|/\rho\le1.0102\,(|\Lambda_2|+E_k(x))\le2.0204\max\{|\Lambda_2|,E_k(x)\}$$
and Lemma~\ref{lem:Lambda}(b) gives the claim. For $L_3$ the
argument is the same, with the additional factor $(1-\varphi^{-4x})/(1-\alpha^{-2x})$, which
differs from $1$ by less than $4x\delta\varphi^{-8}$ and is absorbed in $E_k(x)$.

(b) If $n\le k$, then $P_{n+1}=F_{2n+1}$ by \eqref{eq:kilic}, and by \eqref{eq:zeta},
$\Lambda_1+1=(L_1+1)(1+\zeta_m)$. Hence $$|L_1|\le(|\Lambda_1|+|\zeta_m|)/(1-|\zeta_m|)\le
2.0204\max\{2\cdot5^{-x},|\zeta_m|\}$$ which gives the claim.
\end{proof}

For $L_1$ we apply Theorem~\ref{thm:matveev} in $\Q(\sqrt5)$ with
$(\gamma_1,\gamma_2,\gamma_3)=(2c,\varphi,F_{2n+1})$, $(b_1,b_2,b_3)=(1,2m,-x)$, $B=2m$,
$(A_1,A_2,A_3)=(\log(5+\sqrt5),\log\varphi,2\log F_{2n+1})$ 
and we get
$$\log|L_1|>-4n\log\varphi\cdot K_1'(1+\log2m)$$ with $K_1':=C(2,3)\log(5+\sqrt5)\log\varphi<
9.236\times10^{11}$.

\begin{lemma}\label{lem:Esmall}
For every $k\ge258$ we have $2.022\,E_k(X(k))<10^{-55}$.
\end{lemma}
\begin{proof}
Let $k\ge258$ and $\lambda:=\log k$. We have $C(k,3)\,k^2=c\,k^4(1+\lambda)$ with
$c:=1.4\cdot30^6\cdot3^{4.5}$, and $2c\log\varphi<1.379\cdot10^{11}$. Moreover
$H(k)\le\lambda+1.93$, $\log2+H(k)\le\lambda+2.627$ and
$(1+\lambda)(\lambda+2.627)(\lambda+1.93)\le k^{1.1}$ (at $k=258$ the left side is $<401.2$, the
right side is $>449$, and $k^{1.1}$ divided by the left side increases with $k$). Since
$2\log\varphi<H(k)$, we get
\[
E_1(k)<E_2(k)\le M:=1.379\cdot10^{11}\,k^{5.1},\qquad M\ge2.7\cdot10^{23}.
\]
Put $N(x):=1+\bigl(\log(94.64x)+M(1+\log2x)\bigr)/\log(1+\sqrt2)$, and let $x\ge2$ satisfy one
of the two inequalities of Proposition~5.3. Since $2\log(1+\sqrt2)>\log5$, $N_s(x)\le N(x)$ and
$E_2(k)(1+\log2x)\le N(x)$, in both cases
\[
x\log5<\log2.6+M\,N(x)\bigl(1+\log(N(x)\,x+1)\bigr).\tag{$*$}
\]
Now let $x\ge M^3$ and $u:=\log x$, so that $u\ge160$ and $\log M\le u/3$. Then
$N(x)\le1.135\,M(u+4.7)$ and, since $\log(u+4.7)\le0.035\,u$, $1+\log(N(x)\,x+1)\le1.377\,u$.
Hence the right-hand side of $(*)$ is at most $1.62\,M^2u^2\le1.62\,u^2e^{2u/3}<e^u\log5=x\log5$,
so $(*)$ fails. Therefore $$X(k)\le M^3<2.63\cdot10^{33}k^{15.3}$$ and
\[
2.022\,E_k(X(k))\le G(k):=5.38\cdot10^{33}\,k^{15.3}(0.45k+2.02)\,\varphi^{-2k}.
\]
Since $15.3/k+0.45/(0.45k+2.02)<2\log\varphi$, the function $G$ is decreasing for $k\ge258$, and
$G(309)<10^{-55}$. This proves the claim for $k\ge309$. For $258\le k\le308$ we computed $X(k)$
directly. The largest value of $\log_{10}\bigl(2.022\,E_k(X(k))\bigr)$ is $-55.40$, attained
at $k=258$.
\end{proof}

\begin{proposition}\label{prop:Kabs}
Equation \eqref{eq:main} has no solution satisfying \eqref{eq:remaining} with
$k\ge K_{\rm abs}:=1.0427\times10^{45}$.
\end{proposition}
\begin{proof}
Let $k\ge K_{\rm abs}$ and $x<X:=X(k)$; at $k=K_{\rm abs}$ we have $X<2.99\times10^{398}$. By
Lemma~\ref{lem:Esmall}, Lemma~\ref{lem:transfer} applies. Comparing with Matveev's bound for $L_2$,
the third term $\lambda_*$ cannot be the minimum in $\lambda_2$, because
$$2k\log\varphi>K_2(1+\log3X)+\log\big(2.022\cdot1.011X(0.45k+2.02)\big).$$ Hence, if \eqref{eq:star}
holds, either $2x\log\alpha-\log5.26<K_2(1+\log3x)$, which gives $x<X_0'=1.7634\times10^{13}$
as in Section~\ref{sec:fib}, or
$$n<N_b:=\big(K_2(1+\log3X)+\log(5.26X)\big)/\log\alpha>9.949\times10^{14}$$ if \eqref{eq:star}
fails, $n<1+\log(36.4X)/\log\alpha<N_b$.

\emph{Case $x<X_0'$.} Then Matveev's bound for $L_3$ gives, again because $\lambda_*$ cannot be the
minimum, $n<N_a:=1.32\times10^{29}$. Hence $n\le k$ and $m\le nx+1<N_aX_0'+1<k+1$, so that all
the terms in \eqref{eq:main} are odd-indexed Fibonacci numbers by \eqref{eq:kilic}, and
\eqref{eq:main} becomes \eqref{eq:fib}, which has no solution by Theorem~\ref{thm:fib}.

\emph{Case $n<N_b$.} Then $n\le k$ and $m\le N_bX+1$, and $|\zeta_m|$ is extremely small. Matveev's
bound for $L_1$ and Lemma~\ref{lem:transfer}(b) give $x<X_G=1.048\times10^{30}$. Then $m\le N_bX_G+1\le K_{\rm abs}\le k$. This is the
condition which determines $K_{\rm abs}$. We conclude with Theorem~\ref{thm:fib} as before.

All the inequalities used are of the form $k>F(\log k)$, where $F$ is a polynomial with positive
coefficients.  Since they hold at $k=K_{\rm abs}$, they hold for all
$k\ge K_{\rm abs}$. We also checked them numerically for $k=K_{\rm abs}\cdot10^{t}$,
$t=0,5,10,\ldots,180$.
\end{proof}

We now lower $K_{\rm abs}$ by reduction.
Let $K$ be an upper bound for $k$, with $K\ge258$, and consider $258\le k\le K$; then $x<X(K)$ and,
by Lemma~\ref{lem:Esmall}, Lemma~\ref{lem:transfer} applies. One round of reduction consists of
the following steps, with $\log\alpha\ge2\log\varphi+\log(1-\varphi^{-30})$ and $x\le M$:
\begin{enumerate}[label=\rm(R\arabic*)]
\item Lemma~\ref{lem:dp} applied to $L_2$, with $\tau=\log\sqrt5/\log\varphi$,
$\mu=\log2/\log\varphi$, $u=x-1\le M$, $A=2/\log\varphi$, $B=e$, gives $\lambda_2<W$. Hence
either $x\le X_a:=(W+\log5.26)/(2\log\alpha)$, or $n<(W+\log(5.26M))/\log\alpha$, or
$$k<K_c:=\big(W+\log(2.022\cdot1.011M(0.45K+2.02))\big)/(2\log\varphi)$$ or \eqref{eq:star}
fails.
\item For each $2\le x\le X_a$ and $x-1\le j\le2x-1$ the number $L_3$ is explicit and non-zero;
with $\ell:=-\log\min_j|L_3|$, Lemma~\ref{lem:transfer}(a) gives either
$n<(\ell+\log(3.04x))/\log\alpha$, or $$2k\log\varphi<\ell+\log(2.022\cdot1.011x(0.45K+2.02)).$$
\item Steps (R1)--(R2) give $n\le N$ unless $k$ is smaller than an explicit bound. For each
$2\le n\le N$, Lemma~\ref{lem:dp} applied to $L_1$, with $\tau=2\log\varphi/\log F_{2n+1}$,
$\mu=\log(2c)/\log F_{2n+1}$, $u=m\le NM+1$, $A=2/\log F_{2n+1}$, $B=e$, gives $x\le X_G$ unless
$k$ is smaller than an explicit bound.
\end{enumerate}
So, if $k$ is larger than a bound on $K_1$, let us say the maximum of the bounds met in (R1)--(R3) and of $N$,
we get $n\le N$ and $x\le\max\{X_a,X_G\}$. We repeat (R1)--(R3) with this new bound for $x$ until it
stabilizes, and we get final bounds $n\le N'$, $x\le X'$, valid for $k>K_2$. Finally, if
$k\ge N'X'+1$, then $n\le k$ and $m\le N'X'+1\le k+1$, so that \eqref{eq:main} becomes
\eqref{eq:fib}, which has no solution. Thus a round gives the new bound
$k\le\max\{K_1,K_2,N'X'+1\}$.

\begin{table}[ht]
\centering
\caption{Rounds of reduction for large $k$. Here $K$ is the bound for $k$ at the start of the round. For the first pass of each round we give $\lfloor\log_{10}X(K)\rfloor$, the value $W$ of (R1), $X_a$, the bound $N$ for $n$, and the bound $X_G$ of (R3).
The last pass of every round gives $W=6.40$, $N'=10$, $X'=6$, hence $N'X'+1=61$.}\label{tab:rounds}
\begin{tabular}{@{}rrrrrrrr@{}}
\toprule
round & $K$ & $\log_{10}X(K)$ & $W$ & $X_a$ & $N$ & $X_G$ & new $K$\\
\midrule
1 & $1.0427\times10^{45}$ & $398$ & $923.69$ & $480$ & $1914$ & $580$ & $2021$\\
2 & $2022$ & $57$ & $138.32$ & $72$ & $283$ & $92$ & $297$\\
3 & $298$ & $50$ & $122.13$ & $64$ & $249$ & $80$ & $261$\\
4 & $262$ & $50$ & $120.60$ & $63$ & $246$ & $79$ & $258$\\
5 & $259$ & $50$ & $120.56$ & $63$ & $246$ & $79$ & $257$\\
6 & $258$ & $50$ & $120.55$ & $63$ & $246$ & $79$ & $257$\\
\bottomrule
\end{tabular}
\end{table}

\begin{proposition}\label{prop:rounds}
Equation \eqref{eq:main} has no solution satisfying \eqref{eq:remaining} with $k\ge258$.
\end{proposition}
\begin{proof}
By Proposition~\ref{prop:Kabs} we may assume $k<K_{\rm abs}$. Table~\ref{tab:rounds} shows the
rounds. Each round starts from the bound given by the previous one plus one, and the last
round gives again $k\le257$.
\end{proof}

\section{The cases $2\le k\le257$}\label{sec:smallk}

Now, it  remains is a finite, but very long, computation. We have $256$ values of $k$, each with its own
$\alpha$, its own $g$ and a bound $X(k)$ up to $10^{50}$. The steps are those of
Section~\ref{sec:fib}, in $\Q(\alpha)$ instead of $\Q(\sqrt5)$.

Let $2\le k\le257$. We start with $x<X(k)\le1.08\times10^{50}$ (Proposition~\ref{prop:X}) and
make three passes of the following three steps, with $x\le M$.
\begin{enumerate}[label=\rm(S\arabic*)]
\item With $\Gamma_2=\log2+(1-x)\log g-j\log\alpha$, Lemma~\ref{lem:Lambda}(b) gives, under
\eqref{eq:star},
$$
0<|x\tau-j+\mu|<\frac{5.2}{\log\alpha}e^{-w},\qquad w=\min\{2x\log\alpha,\ n\log\alpha-\log x\},
$$
with $\tau=-\log g/\log\alpha$ and $\mu=\log(2g)/\log\alpha$.
Lemma~\ref{lem:dp} with $u=x\le M$ gives $w<W$. Hence either $x<W/(2\log\alpha)$, or
$n<(W+\log M)/\log\alpha$, or \eqref{eq:star} fails and $n<1+\log(36.4M)/\log\alpha$.
\item For each $2\le x<W/(2\log\alpha)$ and $x-1\le j\le2x-1$ we computed $\Lambda_3$ (it is
non-zero, except for $(k,x,j)=(2,2,2)$, see Lemma~\ref{lem:zero}) and
$n<\log(1.5x/|\Lambda_3|)/\log\alpha$ from Lemma~\ref{lem:Lambda}(b). Together with (S1) this
gives $n\le N$.
\item For each $2\le n\le N$, Lemma~\ref{lem:dp} applied to $\Gamma_1/\log P_{n+1}$, where
$\Gamma_1=\log(2g)+m\log\alpha-x\log P_{n+1}$, with $\tau=\log\alpha/\log P_{n+1}$,
$\mu=\log(2g)/\log P_{n+1}$, $u=m\le NM+1$, $A=4/\log P_{n+1}$ and $B=5$, gives a new bound for
$x$.
\end{enumerate}

\begin{table}[ht]
\centering
\caption{Maximum over $2\le k\le257$ of the quantities met in the three passes.}\label{tab:smallk}
\begin{tabular}{@{}rrrrrr@{}}
\toprule
pass & $M$ & $W$ in (S1) & $x$ in (S2) & $N$ & $x$ from (S3)\\
\midrule
1 & $1.08\times10^{50}$ & $124.45$ & $64$ & $248$ & $82$\\
2 & $82$ & $14.99$ & $7$ & $19$ & $11$\\
3 & $11$ & $9.63$ & $5$ & $12$ & $8$\\
\bottomrule
\end{tabular}
\end{table}

\begin{proposition}\label{prop:smallk}
Let $2\le k\le257$. Every solution of \eqref{eq:main} with $n\ge2$ and $x\ge2$ satisfies $n\le12$
and $x\le8$. For $k\ge3$ there is no such solution, and for $k=2$ all of them have $x=2$ and
$m=2n$.
\end{proposition}
\begin{proof}
The bounds follow from Table~\ref{tab:smallk}.
For each $k$, a search in
exact integer arithmetic over $2\le n\le N_k$, $2\le x\le X_k$, where $N_k\le12$ and $X_k\le8$ are the final
bounds of the third pass for this $k$, comparing $(P_{n+1}^x-P_{n-1}^x)/2$ with the list of the
$P_m$, $m\le N_kX_k+3$, gives the result.
\end{proof}

\section{Proof of Theorem~\ref{thm:main}}\label{sec:proof}

Let $(k,n,m,x)$ be a solution of \eqref{eq:main}. If $x=1$ or $n=1$, or if $k=2$ and $x=2$, we
conclude with Lemmas~\ref{lem:x1}, \ref{lem:n1} and \ref{lem:k2x2}. Otherwise
\eqref{eq:remaining} holds, that is, $n\ge2$, $x\ge2$ and $(k,x)\ne(2,2)$. By
Proposition~\ref{prop:rounds} we have $k\le257$. By Proposition~\ref{prop:smallk}, every solution
with $2\le k\le257$, $n\ge2$ and $x\ge2$ has $k=2$ and $x=2$, which is excluded by
\eqref{eq:remaining}. Hence no solution satisfies \eqref{eq:remaining}, and this proves
Theorem~\ref{thm:main}.\qed

\appendix
\section{Computational details}\label{app:comp}

All computations were done in Python~3. Real numbers were handled with the \texttt{decimal}
module; the precision was chosen in each step as twice the number of digits of the bound $M$
plus $60$ to $120$ digits, so that the products $\tau q$ and $\mu q$ in Lemma~\ref{lem:dp}, where
$q$ is about $6M\cdot10^{20}$, are correct. The dominant root $\alpha$ was computed by Newton's
method on the equation $x^2-3x+1+x^{1-k}=0$ from \eqref{eq:identity}, which is well conditioned
for large $k$, starting from $\varphi^2$, and $g$ from \eqref{eq:g}. The $k$-Pell numbers were
computed exactly. In all the applications of Lemma~\ref{lem:dp} the first convergent with $q>6M$
or the next one gave $\varepsilon>0$.

As a consistency check, independent of the proof, we also verified by a direct search that, for
$3\le k\le100$, $1\le n\le150$, $1\le x\le60$, and for $2\le k\le60$, $1\le n\le80$, $1\le x\le30$,
the only solutions of \eqref{eq:main} are those of Theorem~\ref{thm:main}. The minimum in
Lemma~\ref{lem:zero}(b) was computed with $150$ significant digits.
%

For $k\ge30$ the minimum of $|L_3|$ over $2\le x\le480$ and $x-1\le j\le 2x-1$ is greater than
$8.43\times10^{-4}$, attained at $(x,j)=(457,610)$; so (R2) never meets a small value.

\section*{Declarations}

\paragraph{\textbf{Data availability}} The Python scripts used in this paper are available from
the second  author on request.


\paragraph{\textbf{Use of AI tools}} Claude (Anthropic) was used to find the error in an earlier
version of this work, to
write and run the Python code described in Appendix~\ref{app:comp}. The
authors take full responsibility for the content.

\end{document}